\documentclass[amstex,12pt,russian,amssymb]{article}
\usepackage{mathtext}
\usepackage[cp1251]{inputenc}
\usepackage[T2A]{fontenc}
\usepackage[russian]{babel}
\usepackage[dvips]{graphicx}
\usepackage{amsmath}
\usepackage{amssymb}
\usepackage{amsxtra}
\usepackage{latexsym}
\usepackage{ifthen}

\begin{document}

\newcounter{lemma}
\newcommand{\lemma}{\par \refstepcounter{lemma}%
{\bf Лемма \arabic{lemma}.}}

\newcounter{corollary}
\newcommand{\corollary}{\par \refstepcounter{corollary}%
{\bf Следствие \arabic{corollary}.}}

\newcounter{remark}
\newcommand{\remark}{\par \refstepcounter{remark}%
{\bf Замечание \arabic{remark}.}}

\newcounter{theorem}
\newcommand{\theorem}{\par \refstepcounter{theorem}%
{\bf Теорема \arabic{theorem}.}}

\newcounter{proposition}
\newcommand{\proposition}{\par \refstepcounter{proposition}%
{\bf Предложение \arabic{proposition}.}}

\newcounter{definition}
\newcommand{\definition}{\par \refstepcounter{definition}%
{\bf Определение \arabic{definition}.}}

\renewcommand{\refname}{\centerline{\bf Список литературы}}

\newcommand{\proof}{{\it Доказательство.\,\,}}

\noindent УДК 517.5

{\bf Е.А.~Севостьянов} (Житомирский государственный университет им.\
И.~Франко)

\medskip
{\bf Є.О.~Севостьянов} (Житомирський державний університет ім.\
І.~Франко)

\medskip
{\bf E.A.~Sevost'yanov} (Zhitomir State University of I.~Franko)

\medskip
{\bf О нульмерности предела последовательности отображений,
удовлетворяющих одному модульному неравенству}

{\bf Про нульвимірність границі послідовності відображень, що
задовольняють одну модульну нерівність }

{\bf On the lightness of the limit of sequence of mappings
satisfying some modular inequality}

\medskip
Статья посвящена изучению свойств одного класса
прос\-тра\-н\-с\-т\-вен\-ных отображений, более общих, чем
отображения с ограниченным искажением. Показано, что локально
равномерный предел последовательности отображений $f:D\rightarrow
{\Bbb R}^n$ области $D\subset{\Bbb R}^n,$ $n\geqslant 2,$
удовлетворяющих одному неравенству относительно $p$-модуля семейств
кривых, является нульмерным. Указанное утверждение обобщает
известную теорему об открытости и дискретности равномерного предела
последовательности отображений с ограниченным искажением.

\medskip
Статтю присвячено вивченню властивостей одного класу просторових
відображень більш загальних, ніж відображення з обмеженим
спотворенням. Доведено, що локально рівномірною границею
послідовності відображень $f:D\rightarrow {\Bbb R}^n$ області
$D\subset{\Bbb R}^n,$ $n\geqslant 2,$ які задовольняють одну
нерівність відносно $p$-модуля сімей кривих, є нульвимірним. Вказане
твердження узагальнює відому теорему про відкритість і дискретність
рівномірної границі послідовності відображень з обмеженим
спотворенням.

\medskip
A paper is devoted to study of one class of space mappings which are
more general than mappings with bounded distortion. It is showed
that a locally uniformly limit of a sequence of mappings
$f:D\rightarrow {\Bbb R}^n$ of domain $D\subset{\Bbb R}^n,$
$n\geqslant 2,$ satisfying one inequality with respect to
$p$-modulus of families of curves, is light. The above statement is
a generalization of well-known theorem about openness and
discreteness of uniformly limit of a sequence of mappings with
bounded distortion.

\newpage
{\bf 1. Введение.} Настоящая заметка посвящена изучению отображений
с ограниченным и конечным искажением, активно изучаемых в последнее
время (см., напр., \cite{ARS}--\cite{Va}). Здесь же хотелось бы
указать на работы, связанные с классами Соболева и пространствами
Карно--Каратедори новосибирской школы математиков
(см.~\cite{Vo$_1$}--\cite{UV}).

Напомним некоторые определения. Всюду далее $D$ -- область в ${\Bbb
R}^n,$ $n\ge 2,$ $m$ -- мера Лебега ${\Bbb R}^n,$ запись
$f:D\rightarrow {\Bbb R}^n$ предполагает, что отображение $f,$
заданное в области $D,$ непрерывно. Как обычно, мы пишем $f\in
W^{1,n}_{loc}(D),$ если все координатные функции
$f=(f_1,\ldots,f_n)$ обладают обобщёнными частными производными
первого порядка, которые локально интегрируемы в $D$ в степени $n.$
Отображение $f:D\rightarrow {\Bbb R}^n$ называется {\it дискретным},
если прообраз $f^{-1}\left(y\right)$ каждой точки $y\,\in\,{\Bbb
R}^n$ состоит из изолированных точек и {\it открытым}, если образ
любого открытого  множества $U\subset D$ является открытым
множеством в ${\Bbb R}^n.$ Напомним, что отображение $f:D\rightarrow
{\Bbb R}^n$ называется {\it отображением с ограниченным искажением},
если выполнены следующие условия:

\noindent 1) $f\in W_{loc}^{1,n},$ 2) якобиан $J(x,f)$ отображения
$f$ в точке $x\in D$ сохраняет знак почти всюду в $D,$ 3) $\Vert
f^{\,\prime}(x) \Vert^n \le K \cdot |J(x,f)|$ при почти всех $x\in
D$ и некоторой постоянной $K<\infty,$ где, как обычно, $\Vert
f^{\,\prime}(x)\Vert:=\sup\limits_{h\in {\Bbb R}^n:
|h|=1}|f^{\,\prime}(x)h|,$
см., напр., \cite[$\S\, 3,$ гл. I]{Re}, либо \cite[определение 2.1,
разд.~2, гл.~I]{Ri}. Начало интенсивных исследований
пространственных отображений с ограниченным искажением положено Ю.
Г. Решетняком. В его работах, в частности, доказаны открытость и
дискретность отображений $f$ с ограниченным искажением, см.
\cite[теоремы~6.3 и 6.4, $\S\, 6,$ гл.~II]{Re}.

\medskip
Обратим внимание на следующий результат, также принадлежащий
Решетняку: {\it если последовательность $f_m$ отображений с
ограниченным искажением, имеющая общую постоянную квазиконформности
$K\geqslant 1,$ сходится локально равномерно к отображению $f,$ то
$f$ является отображением с ограниченным искажением, в частности,
$f$ открыто и дискретно} (см. \cite{Re}). В настоящей статье
указанный результат в несколько модифицированном варианте
доказывается нами для более широкого класса отображений. Для
формулировки этого результата рассмотрим ряд дополнительных
определений, связанных с геометрическим истолкованием отображений с
ограниченным искажением.

\medskip
Здесь и далее {\it кривой} $\gamma$ мы называем непрерывное
отображение отрезка $[a,b]$ (либо открытого интервала $(a,b),$ а
также полуоткрытых интервалов вида $[a,b),$ $(a,b]$) в ${\Bbb R}^n,$
$\gamma:[a,b]\rightarrow {\Bbb R}^n.$ Под семейством кривых $\Gamma$
подразумевается некоторый фиксированный набор кривых $\gamma,$ а
$f(\Gamma)=\left\{f\circ\gamma|\gamma\in\Gamma\right\}.$ Следующие
определения могут быть найдены, напр., в \cite[разд.~1--6]{Va}.
Борелева функция $\rho:{\Bbb R}^n\,\rightarrow [0,\infty]$
называется {\it допустимой} для семейства $\Gamma$ кривых $\gamma$ в
${\Bbb R}^n,$ если криволинейный интеграл первого рода
$\int\limits_{\gamma}\rho (x) |dx|$ удовлетворяет условию
$\int\limits_{\gamma}\rho (x) |dx| \geqslant 1$ для всех кривых $
\gamma \in \Gamma.$ В этом случае мы пишем: $\rho \in {\rm adm}
\,\Gamma.$ Пусть $p\geqslant 1,$ тогда {\it $p$ -- модулем}
семейства кривых $\Gamma $ называется величина
$M_p(\Gamma)=\inf\limits_{\rho \in \,{\rm adm}\,\Gamma}
\int\limits_{{\Bbb R}^n} \rho^p(x)dm(x).$ Говорят, что семейство
кривых $\Gamma_1$ {\it минорируется} семейством $\Gamma_2,$ пишем
$\Gamma_1\,>\,\Gamma_2,$
если для каждой кривой $\gamma\,\in\,\Gamma_1$ существует подкривая,
которая принадлежит семейству $\Gamma_2.$
В этом случае,
\begin{equation}\label{eq32*A}
\Gamma_1
> \Gamma_2 \quad \Rightarrow \quad M_p(\Gamma_1)\le M_p(\Gamma_2)
\end{equation} (см. \cite[теорема~6.4, гл.~I]{Va}).

\medskip
В относительно недавней работе \cite{Sev$_1$} установлена открытость
и дискретность отображений $f$ области $D\subset{\Bbb R}^n,$
$n\geqslant 2,$ в ${\Bbb R}^n,$ удовлетворяющих для произвольной
функции $\rho_*\in {\rm adm}\,f(\Gamma)$ оценке вида
\begin{equation} \label{eq2*B}
M(\Gamma )\leqslant \int\limits_{f(D)} Q(y)\cdot \rho_*^n (y) dm(y)
\end{equation}
относительно конформного модуля семейств кривых
$M(\Gamma):=M_n(\Gamma)$ и некоторой заданной функции $Q:{\Bbb
R}^n\rightarrow [0, \infty],$ $Q(x)\equiv 0$ при всех $x\in {\Bbb
R}^n\setminus f(D).$
В настоящей работе мы покажем, что локально равномерным пределом
отображений, удовлетворяющих оценке вида (\ref{eq2*B}), может быть
только нульмерное или постоянное отображение. Заметим, при этом, что
отображения с ограниченным искажением $f:D\rightarrow {\Bbb R}^n$
удовлетворяют неравенству (\ref{eq2*B}), поскольку для них всегда
\begin{equation}\label{eq2}
M(\Gamma)\leqslant N(f, A)\cdot K\cdot M(f(\Gamma))\,,
\end{equation}
где $$N(y, f, A)\,=\,{\rm card}\,\left\{x\in A\,|\,
f(x)=y\right\}\,, \qquad N(f, A)\,=\,\sup\limits_{y\in{\Bbb
R}^n}\,N(y, f, A)\,,$$
$A$ -- произвольное борелевское множество в $D,$ а $K\geqslant 1$ --
некоторая постоянная, которая может быть вычислена как $K={\rm ess
\sup}\, K_O(x, f),$ $K_O(x, f)=\Vert f^{\,\prime}(x)\Vert^n/J(x, f)$
при $J(x, f)\ne 0,$ $K_O(x, f)=1$ при $f^{\,\prime}(x)=0,$ и $K_O(x,
f)=\infty$ при $f^{\,\prime}(x)\ne 0,$ но $J(x, f)=0$ (см.
\cite[теорема~6.7, гл.~II]{Ri}). Таким образом, основной результат
данной заметки, сформулированный ниже, есть некая модификация
результата Решетняка о пределе последовательности отображений с
ограниченным искажением, соответствующая случаю, когда $Q(y)$
ограничено в (\ref{eq2*B}).

\medskip
Теперь приведём ещё неоторые вспомогательные сведения. Множество
$H\subset \overline{{\Bbb R}^n},$ $\overline{{\Bbb R}^n}:={\Bbb
R}^n\{\infty\},$ будем называть {\it всюду разрывным}, если любая
его компонента связности вырождается в точку; в этом случае пишем
${\rm dim\,}H=0,$ где ${\rm dim}$ обозначает {\it топологическую
размерность} множества $H,$ см. \cite{HW}. Отображение
$f:D\rightarrow \overline{{\Bbb R}^n}$ называется {\it нульмерным},
если ${\rm dim\,}\{f^{\,-1}(y)\}=0$ для каждого $y\in
\overline{{\Bbb R}^n}.$

Пусть $Q:D\rightarrow [0,\infty]$ -- измеримая по Лебегу функция,
тогда $q_{x_0}(r)$ означает среднее интегральное значение $Q(x)$ над
сферой $S(x_0, r),$
\begin{equation}\label{eq32*}
q_{x_0}(r):=\frac{1}{\omega_{n-1}r^{n-1}}\int\limits_{|x-x_0|=r}Q(x)\,dS\,,
\end{equation}
где $dS$ -- элемент площади поверхности $S.$ Будем говорить, что
функция ${\varphi}:D\rightarrow{\Bbb R}$ имеет {\it конечное среднее
колебание} в точке $x_0\in D$, пишем $\varphi\in FMO(x_0),$ если
%\begin{equation} \label{eq20*}
%
$$\overline{\lim\limits_{\varepsilon\rightarrow 0}}\ \
\frac{1}{\Omega_n\varepsilon^n} \ \ \int\limits_{B(
x_0,\,\varepsilon)}
|{\varphi}(x)-\overline{{\varphi}}_{\varepsilon}|\ dm(x)\, <\,
\infty,
$$
%\end{equation}
%
где
%
%\begin{equation}\label{eq1.1.23}
$\overline{{\varphi}}_{\varepsilon}\,=\,\frac{1}{\Omega_n\varepsilon^n}\int\limits_{B(
x_0,\,\varepsilon)} {\varphi}(x)\ dm(x)$ (см., напр., \cite[разд.
6.1]{MRSY}). Заметим, что все ограниченные функции $\varphi$ --
конечного среднего колебания. Для точки $y_0\in {\Bbb R}^n$ и чисел
$0<r_1<r_2<\infty$ обозначим
\begin{equation}\label{eq1**}
A(y_0, r_1,r_2)=\left\{ y\,\in\,{\Bbb R}^n:
r_1<|y-y_0|<r_2\right\}\,.\end{equation}
Пусть $E,$ $F\subset\overline{{\Bbb R}^n}$ -- произвольные
множества. Обозначим через $\Gamma(E,F,D)$ семейство всех кривых
$\gamma:[a,b]\rightarrow \overline{{\Bbb R}^n},$ которые соединяют
$E$ и $F$ в $D\,,$ т.е. $\gamma(a)\in E,\gamma(b)\in\,F$ и
$\gamma(t)\in D$ при $t \in (a, b).$ Если $f:D\rightarrow {\Bbb
R}^n$ -- заданное отображение, то для фиксированного $y_0\in f(D)$ и
произвольных $0<r_1<r_2<\infty$ обозначим через $\Gamma_f(y_0, r_1,
r_2)$ семейство всех кривых $\gamma$ в области $D$ таких, что
$f(\gamma)\in \Gamma(S(y_0, r_1), S(y_0, r_2), A(y_0,r_1,r_2)).$
Рассмотрим вместо (\ref{eq2*B}) неравенство
\begin{equation} \label{eq2*A}
M_p(\Gamma_f(y_0, r_1, r_2))\leqslant \int\limits_{f(D)} Q(y)\cdot
\eta^p (y) dm(y)
\end{equation}
выполненное для любой измеримой по Лебегу функции $\eta:
(r_1,r_2)\rightarrow [0,\infty ]$ такой, что
\begin{equation}\label{eqA2}
\int\limits_{r_1}^{r_2}\eta(r) dr\geqslant 1\,.
\end{equation}
Отметим, что даже при $p=n$ неравенство (\ref{eq2*A}) является более
слабым, чем (\ref{eq2*B}). Действительно, если для произвольной
функции $\rho_*\in {\rm adm}\,f(\Gamma)$ имеет место неравенство
(\ref{eq2*B}), то возьмём произвольную измеримую по Лебегу функцию
$\eta: (r_1,r_2)\rightarrow [0,\infty ],$ удовлетворяющую
соотношению (\ref{eqA2}) и определим функцию
$\rho_*(y):=\eta(|y-y_0|).$ Эта функция допустима для семейства
кривых $\Gamma(S_1, S_2, A),$ соединяющих сферы $S_1=S(y_0, r_1)$ и
$S_2=S(y_0, r_2),$ поскольку согласно \cite[теорема~5.7]{Va}
интеграл от произвольной радиальной функции $\Psi(|y-y_0|)$ по
(локально спрямляемой) кривой, соединяющей сферы $S(y_0, r_1)$ и
$S(y_0, r_2)$ не меньше, чем соответствующий интеграл по отрезку
$(r_1, r_2)$ от функции $\Psi(t)$ и, значит,
$\int\limits_{\gamma}\,\rho_{*}(y)\,|dy|\,\geqslant\int\limits_{r_1}
^{r_2}\,\eta(t)dt\geqslant 1$
для произвольной кривой $\gamma\in \Gamma(S_1, S_2, A).$ Значит,
определённую выше функцию $\rho_*$ можно подставить в соотношение
(\ref{eq2*B}), откуда и последует справедливость неравенства
(\ref{eq2*A}).

\medskip
Основной результат настоящей статьи заключает в себе следующая

\medskip
\begin{theorem}\label{th3}{\sl\,Пусть $p\in [n-1, n],$ $Q:{\Bbb R}^n\rightarrow (0, \infty)$ --
измеримая по Лебегу функция, $f_m:D\,\rightarrow\,{\Bbb R}^n,$
$n\geqslant 2$ -- последовательность отображений, удовлетворяющих
(\ref{eq2*A})--(\ref{eqA2}), и сходящаяся локально равномерно к
некоторому отображению $f:D\rightarrow {\Bbb R}^n.$ Пусть функция
$Q(y),$ кроме того, удовлетворяет хотя бы одному из следующих
условий:

1) $Q\in FMO(y_0)$ в произвольной точке $y_0\in f(D),$

2) $q_{y_0}(r)\,=\,O\left(\left[\log{\frac1r}\right]^{n-1}\right)$
при $r\rightarrow 0$ и при всех $y_0\in f(D),$ где функция
$q_{y_0}(r)$ определена равенством (\ref{eq32*}),

3) для каждого $y_0\in f(D)$ найдётся некоторое число
$\delta(y_0)>0,$ такое что при достаточно малых $\varepsilon>0$
\begin{equation}\label{eq5**}
\int\limits_{\varepsilon}^{\delta(y_0)}\frac{dt}{t^{\frac{n-1}{p-1}}q_{y_0}^{\frac{1}{p-1}}(t)}<\infty,
\qquad
\int\limits_{0}^{\delta(y_0)}\frac{dt}{t^{\frac{n-1}{p-1}}q_{y_0}^{\frac{1}{p-1}}(t)}=\infty\,.
\end{equation} Тогда отображение $f$ либо нульмерно,
либо постоянно в $D.$

Если, кроме того, отображение $f$ сохраняет ориентацию, то $f$
открыто и дискретно.}
\end{theorem}

\medskip
{\bf 2. Формулировка и доказательство основной леммы.} Для
проведения доказательств основных результатов нам необходимы
не\-которые сведения из теории общих метрических пространств.
Напомним, что связный компакт $C$ метрического пространства $X$
называется {\it континуумом}. Пусть $(X, \mu)$ -- метрическое
пространство с мерой $\mu.$ Определим {\it функцию Лёвнера
$\phi_n:(0, \infty)\rightarrow [0, \infty)$ на $X$} по следующему
правилу:
$$\phi_n(t)=\inf\{M_n(\Gamma(E, F, X)): \Delta(E, F)\leqslant t\}\,,$$
где $\inf$ берётся по всем произвольным невырожденным
непересекающимся континуумам $E, F$ в $X,$ относительно которых
величина $\Delta(E, F)$ определяется как
$$\Delta(E, F):=\frac{{\rm dist}\,(E,
F)}{\min\{{\rm diam\,}E, {\rm diam\,}F\}}\,.$$
Пространство $X$ называется {\it пространством Лёвнера,} если
функция $\phi_n(t)$ положительна при всех положительных значениях
$t$ (см. \cite[разд.~2.5]{MRSY} либо \cite[гл.~8]{He}). Заметим, что
пространство ${\Bbb R}^n,$ равно как и единичный шар ${\Bbb
B}^n\subset {\Bbb R}^n,$ являются пространствами Лёвнера
относительно стандартной евклидовой метрики и стандартной лебеговой
меры (см. \cite[теорема~8.2 и пример~8.24(a)]{He}). Заметим, что в
пространствах Лёвнера $X$ условие $\mu(B(x_0, r))\geqslant C\cdot
r^n$ выполняется для каждой точки $x_0\in X,$ некоторой постоянной
$C$ и всех $r<{\rm diam}\,X.$ Пространство $X$ будет называться {\it
геодезическим}, если любые две его точки могут быть соединены
кривой, длина которой равна расстоянию между указанными точками. В
частности, ${\Bbb B}^n$ -- геодезическое пространство. Следующее
определение см., напр., в \cite[разд.~1.4, гл.~I]{He}, либо
\cite[раздел~1]{AS}). Говорят, что метрическое пространство $(X,
\rho)$ с мерой $\mu$ является {\it пространством с условием удвоения
меры,} если существует постоянная $C>0$ такая, что для всех $r>0$ и
всех $x_0\in X$ выполняется следующее условие: $\mu(B(x_0, 2r))\le
C\cdot \mu(B(x_0, r)).$ Легко видеть, что произвольная ограниченная
евклидова область $D$ удовлетворяет условию удвоения меры. Следуя
\cite[раздел 7.22]{He} будем говорить, что борелева функция
$\rho\colon  X\rightarrow [0, \infty]$ является {\it верхним
градиентом} функции $u\colon X\rightarrow {\Bbb R},$ если для всех
спрямляемых кривых $\gamma,$ соединяющих точки $x$ и $y\in X,$
выполняется неравенство $|u(x)-u(y)|\leqslant
\int\limits_{\gamma}\rho\,|dx|,$ где, как обычно,
$\int\limits_{\gamma}\rho\,|dx|$ обозначает линейный интеграл от
функции $\rho$ по кривой $\gamma.$ Будем также говорить, что в
указанном пространстве $X$ выполняется {\it $(1; p)$-неравенство
Пуанкаре,} если найдётся постоянная $C\geqslant 1$ такая, что для
каждого шара $B\subset X,$ произвольной ограниченной непрерывной
функции $u\colon X\rightarrow {\Bbb R}$ и любого её верхнего
градиента $\rho$ выполняется следующее неравенство:
$$\frac{1}{\mu(B)}\int\limits_{B}|u-u_B|d\mu(x)\leqslant C\cdot({\rm diam\,}B)\left(\frac{1}{\mu(B)}
\int\limits_{B}\rho^p d\mu(x)\right)^{1/p}\,.$$
Метрическое пространство $(X, d, \mu)$ назовём {\it $n$-регулярным
по Альфорсу,} если при каждом $x_0\in X,$ некоторой постоянной
$C\geqslant 1$ и всех $R<{\rm diam}\,X$
$$\frac{1}{C}R^{n}\leqslant \mu(B(x_0, R))\leqslant CR^{n}\,.$$

\medskip
Справедливо следующее утверждение.

\medskip
\begin{proposition}\label{pr1}
Единичный шар ${\Bbb B}^n$ является $n$-регулярным по Альфорсу
ме\-трическим пространством, в котором выполнено $(1;
p)$-неравенство Пуанкаре. Более того, для любых двух континуумов $E,
F\subset {\Bbb B}^n$ и произвольного $p\in [n-1, n]$
\begin{equation}\label{eq3}
M_p(\Gamma(E, F, {\Bbb B}^n))>0\,.
\end{equation}
\end{proposition}

\begin{proof}
То, что ${\Bbb B}^n$ является $n$-регулярным по Альфорсу,
непосредственно следует из сделанных выше замечаний. Согласно этим
же замечаниям  пространство ${\Bbb B}^n$ является геодезическим, и
является пространством Лёвнера, поэтому в нём выполняется $(1;
p)$-неравенство Пуанкаре (см. \cite[теоремы~9.8 и 9.5]{He}).
Соотношение (\ref{eq3}), в таком случае, есть результат
\cite[следствие~4.8]{AS}.
\end{proof}$\Box$

\medskip
Следующая лемма включает в себя основной результат настоящей работы
в наиболее общей ситуации.

\medskip
\begin{lemma}\label{lem1}
{\sl\,Пусть $p\in [n-1, n],$ $Q:{\Bbb R}^n\rightarrow (0, \infty)$
-- измеримая по Лебегу функция, $f_m:D\,\rightarrow\,{\Bbb R}^n,$
$n\geqslant 2$ -- последовательность отображений, удовлетворяющих
оценкам (\ref{eq2*A})--(\ref{eqA2}) и сходящаяся локально равномерно
к некоторому отображению $f:D\rightarrow {\Bbb R}^n.$ Далее,
предположим, что для каждого $y_0\in D$ найдётся $\varepsilon_0>0,$
для которого выполнено соотношение
\begin{equation} \label{eq4!}
\int\limits_{A(y_0, \varepsilon,
\varepsilon_0)}Q(y)\cdot\psi^p(|y-y_0|) \
dm(y)\,=\,o\left(I^p(\varepsilon, \varepsilon_0)\right)
\end{equation}
для некоторой борелевской функции $\psi(t):(0,\infty)\rightarrow
[0,\infty],$ такой что
\begin{equation} \label{eq5}
0< I(\varepsilon,
\varepsilon_0):=\int\limits_{\varepsilon}^{\varepsilon_0}\psi(t)dt <
\infty
\end{equation}
при всех $\varepsilon\in(0,\varepsilon_0),$ где
$$A(y_0, \varepsilon, \varepsilon_0)=\{y\in {\Bbb R}^n: \varepsilon<|y-y_0|<\varepsilon_0\}\,.$$
Тогда отображение $f$ либо постоянно, либо нульмерно.

Если, кроме того, отображение $f$ сохраняет ориентацию, то $f$
открыто и дискретно.}
\end{lemma}

\medskip
\begin{remark}\label{rem1}
В условиях леммы \ref{lem1}, можно считать, что для произвольного
фиксированного $A,$ такого что $0<A<\varepsilon_0,$ и всех
$\varepsilon\in (0, A),$ выполняется условие вида
$\int\limits_{\varepsilon}^{A}\,\psi(t) dt>0.$ Действительно, из
того, что $Q(x)>0$ п.в., а также соотношений (\ref{eq4!}) и
(\ref{eq5}) следует, что $\int\limits_{\varepsilon}^{A}\,\psi(t)
dt\rightarrow \infty$ при $\varepsilon\rightarrow 0,$ поскольку
величина интеграла слева в (\ref{eq4!}) увеличивается при уменьшении
$\varepsilon.$
\end{remark}

\medskip
%\begin{proof}
{\it Доказательство леммы \ref{lem1}.} Если $f\equiv const,$
доказывать нечего. Пусть $f\not\equiv const.$ Предположим противное,
а именно, $f$ не является нульмерным отображением. Тогда найдётся
$y_0\in {\Bbb R}^n,$ такое что множество $\{f^{\,-1}(y_0)\}$
содержит невырожденный континуум $C\subset \{f^{\,-1}(y_0)\}.$
Поскольку по предположению $f\not\equiv y_0,$ найдётся точка $a,$
принадлежащая этому континууму, в любой окрестности $U$ которой
имеются точки, ему не принадлежащие. Можно считать, что $U={\Bbb
B}^n$ и $x_0\in {\Bbb B}^n: f(x_0)\ne y_0.$ По теореме о сохранении
знака найдётся $x_0\in {\Bbb B}^n$ и $\delta_0>0:$ $\overline{B(x_0,
\delta_0)}\subset {\Bbb B}^n$ и
\begin{equation}\label{eq7*}
f(x)\ne y_0\qquad %\text{при\,\, всех}
\forall\quad x\in \overline{B(x_0, \delta_0)}\,.
\end{equation}
В силу \cite[лемма~1.15]{Na} при $p=n$ и предложения \ref{pr1} при
$p\in [n-1, n)$ будем иметь
\begin{equation}\label{eq4*}
M_p\left(\Gamma\left(C, \overline{B(x_0, \delta_0)}, {\Bbb
B}^n\right)\right)>0\,.
\end{equation}
Зафиксируем достаточно большое $m\in {\Bbb N}$ и рассмотрим
семейство кривых $$f_m\left(\Gamma\left(C, \overline{B(x_0,
\delta_0)}, {\Bbb B}^n\right)\right)\,.$$ Заметим, что в силу
локально равномерной сходимости $f_m$ к $f$ может быть построена
подпоследовательность $f_{m_k}$ такая, что $|f_{m_k}(x)-y_0|<1/2^k$
при всех $k\in {\Bbb N}$ и всех $x\in C.$ С другой стороны,
$f(\overline{B(x_0, \delta_0)})$ -- компакт в ${\Bbb R}^n,$ поэтому
${\rm dist}\,(y_0, f(\overline{B(x_0, \delta_0)}))\geqslant
\sigma_0>0.$ Поскольку $f_m$ сходится к $f$ локально равномерно,
$$|f_m(x)-y_0|=|f_m(x)-f(x)+f(x)-y_0|\geqslant |f(x)-y_0|-
|f_m(x)-f(x)|\geqslant \sigma_0/2$$ при всех $x\in \overline{B(x_0,
\delta_0)}$ и всех $m\geqslant m_0.$

В таком случае, каждая кривая $\gamma\in f_{m_k}\left(\Gamma\left(C,
\overline{B(x_0, \delta_0)}, {\Bbb B}^n\right)\right)$ имеет
подкривую $\gamma^{\,\prime}\in \Gamma(S(y_0, 1/2^k), S(y_0,
\sigma_0/2), A(y_0, 1/2^k, \sigma_0/2))$ при достаточно больших
$k\geqslant k_0$ (см. \cite[предложение~13.3]{MRSY}). Отсюда
$\Gamma\left(C, \overline{B(x_0, \delta_0)}, {\Bbb
B}^n\right)>\Gamma_{f_{m_k}}(y_0, 1/2^k, \sigma_0/2)$ и, значит, в
силу (\ref{eq32*A})
\begin{equation}\label{eq1}
M_p\left(\Gamma\left(C, \overline{B(x_0, \delta_0)}, {\Bbb
B}^n\right)\right)\leqslant M_p(\Gamma_{f_{m_k}}(y_0, 1/2^k,
\sigma_0/2))\,.
\end{equation}
Рассмотрим следующую функцию
$$\eta_{k}(t)=\left\{
\begin{array}{rr}
\psi(t)/I(1/2^k, \sigma_0/2), & t\in [1/2^k,
\sigma_0/2],\\
0,  & t\in {\Bbb R}\setminus [1/2^k, \sigma_0/2]\,,
\end{array}
\right. $$
где $I(1/2^k, \sigma_0/2)=\int\limits_{1/2^k}^{\sigma_0/2}\,\psi (t)
dt.$ Заметим, что функция $\eta_{k}$ удовлетворяет условию вида
(\ref{eqA2}) при $r_1=1/2^k$ и $r_2=\sigma_0/2.$ Тогда согласно
неравенствам (\ref{eq2*A}) и (\ref{eq1}) мы получаем, что
\begin{equation}\label{eq6}
M_p\left(\Gamma\left(C, \overline{B(x_0, \delta_0)}, {\Bbb
B}^n\right)\right)\leqslant \frac{1}{I^p(1/2^k,
\sigma_0/2)}\int\limits_{A(y_0, 1/2^k,
\sigma_0/2)}Q(y)\psi^p(|y-y_0|)dm(y)\,\leqslant
\end{equation}
$$\leqslant \frac{C}{I^p(1/2^k,
\varepsilon_0)}\int\limits_{A(y_0, 1/2^k,
\sigma_0/2)}Q(y)\psi^p(|y-y_0|)dm(y) \rightarrow 0$$
при $k\rightarrow \infty.$
Однако, соотношение (\ref{eq6}) противоречит неравенству
(\ref{eq4*}). Полученное противоречие говорит о том, что отображение
$f$ нульмерно.

Пусть дополнительно известно, что $f$ сохраняет ориентацию. Тогда
открытость и дискретность отображения $f$ вытекает из
\cite[следствие, с.~333]{TY}~$\Box$

{\bf 3. О доказательстве основного результата.} {\it Доказательство
теоремы \ref{th3}} непосредственно вытекает из леммы \ref{lem1} и
следующих соображений. В случае 1), когда $Q\in FMO(y_0),$
необходимо рассмотреть функцию
$\psi(t)=\left(t\,\log{\frac1t}\right)^{-n/p}>0$ и применить к ней
утверждение леммы \ref{lem1}. Тогда, ввиду \cite[следствие~6.3,
гл.~6]{MRSY}, условие $Q\in FMO(y_0)$ влечёт, что при достаточно
малых $\varepsilon<\varepsilon_0$
\begin{equation}\label{eq31*}
\int\limits_{\varepsilon<|y-y_0|<\varepsilon_0}Q(y)\cdot\psi^p(|y-y_0|)
\ dm(y)\,=\,O \left(\log\log \frac{1}{\varepsilon}\right)\,,\quad
\varepsilon\rightarrow 0\,.
\end{equation}
Заметим, что величина $I(\varepsilon, \varepsilon_0),$ определённая
в лемме \ref{lem1}, может быть оценена следующим образом:
\begin{equation}\label{eqlogest}
I(\varepsilon,
\varepsilon_0)=\int\limits_{\varepsilon}^{\varepsilon_0}\psi(t) dt
>\log{\frac{\log{\frac{1}
{\varepsilon}}}{\log{\frac{1}{\varepsilon_0}}}}.
\end{equation}
В таком случае, из условия (\ref{eq31*}) следует, что
$$ \frac{1}{I^p(\varepsilon,
\varepsilon_0)}\int\limits_{\varepsilon<|y-y_0|<\varepsilon_0}
Q(x)\cdot\psi^p(|y-y_0|)dm(y)\leqslant
C\left(\log\log\frac{1}{\varepsilon}\right)^{1-p}\rightarrow 0,
\quad \varepsilon\rightarrow 0\,,
$$
что и завершает доказательство теоремы в случае 1), поскольку здесь
выполнены условия (\ref{eq4!})--(\ref{eq5}) леммы \ref{lem1}.
Заметим, что случай 2) является частным случаем ситуации 3), поэтому
для завершения доказательства теоремы достаточно рассмотреть только
случай 3). В этом случае полагаем
\begin{equation}\label{eq9}
I=I(\varepsilon,\varepsilon_0)=\int\limits_{\varepsilon}^{\varepsilon_0}\
\frac{dr}{r^{\frac{n-1}{p-1}}q_{y_0}^{\frac{1}{p-1}}(r)}\,.
\end{equation}
Для произвольных $0<\varepsilon<\varepsilon_0<1$ рассмотрим функцию
\begin{equation}\label{eq1*****}
\psi(t)= \left \{\begin{array}{rr}
1/[t^{\frac{n-1}{p-1}}q_{y_0}^{\frac{1}{p-1}}(t)]\ , & \ t\in
(\varepsilon,\varepsilon_0)\ ,
\\ 0\ ,  &  \ t\notin (\varepsilon,\varepsilon_0)\ .
\end{array} \right.
\end{equation}
Заметим, что функция $\psi$ удовлетворяет всем условиям леммы
\ref{lem1}.
По теореме Фубини (см. \cite[теорема~8.1, гл.~III]{Sa}) имеем
$\int\limits_{\varepsilon<|y-y_0|<\varepsilon_0}
Q(y)\cdot\psi^p(|y-y_0|)\, dm(y)=\omega_{n-1}\cdot I(\varepsilon,
\varepsilon_0)$ (где $\omega_{n-1}$ -- площадь единичной сферы
${\Bbb S}^{n-1}$ в ${\Bbb R}^n$). Вывод: выполнены условия
(\ref{eq4!})--(\ref{eq5}) леммы \ref{lem1}, что окончательно и
доказывает теорему. $\Box$

\medskip
{\bf 4. Основные следствия.} Из леммы \ref{lem1} получаем также
следующие утверждения.

\medskip
\begin{corollary}\label{cor1}
{\sl\,В условиях теоремы \ref{th3} отображение $f$ является
нульмерным либо постоянным, если в каждой точке $y_0\in f(D)$ при
некотором $\delta_0=\delta_0(y_0)$ и всех достаточно малых
$\varepsilon>0$ выполнено условие
$\int\limits_{\varepsilon}^{\delta_0}\frac{dt}{tq_{y_0}^{\frac{1}{n-1}}(t)}<\infty$
и, кроме того,
$\int\limits_{0}^{\delta_0}\frac{dt}{tq_{y_0}^{\frac{1}{n-1}}(t)}=\infty\,.$
}
\end{corollary}

\medskip\begin{proof}
Аналогично доказательству пункту 3) теоремы \ref{th3} рассмотрим
функцию $\psi(t)\quad=\quad \left \{\begin{array}{rr}
\left(1/[tq^{\frac{1}{n-1}}_{y_0}(t)]\right)^{n/p}\ , & \ t\in
(\varepsilon, \varepsilon_0)\ ,
\\ 0\ ,  &  \ t\notin (\varepsilon,
\varepsilon_0)\ ,
\end{array} \right.$ Рассуждая также, как в указанном случае,
получаем требуемое утверждение.
\end{proof}$\Box$

\medskip
Отдельно рассмотрим случай $p\in [n-1, n).$

\medskip
\begin{corollary}\label{cor3}
{\sl\, Пусть $p\in [n-1, n),$ тогда в условиях теоремы \ref{th3}
отображение $f$ является нульмерным либо постоянным, если функция
$Q$ удовлетворяет условию $Q\in L_{loc}^s({\Bbb R}^n)$ при некотором
$s\geqslant\frac{n}{n-p}.$ }
\end{corollary}

\medskip
\begin{proof}
Зафиксируем произвольным образом $0<\varepsilon_0<\infty$ и для
произвольного $y_0\in f(D)$ положим $G:=B(y_0, \varepsilon_0)$ и
$\psi(t):=1/t.$ Заметим, что указанная функция $\psi$ удовлетворяет
неравенствам (\ref{eq5}), так что остаётся проверить лишь
справедливость условия (\ref{eq4!}). Применяя неравенство Гёльдера,
будем иметь
$$\int\limits_{\varepsilon<|x-b|<\varepsilon_0}\frac{Q(x)}{|x-b|^p} \
dm(x)\leqslant$$
\begin{equation}\label{eq13}
\leqslant
\left(\int\limits_{\varepsilon<|x-b|<\varepsilon_0}\frac{1}{|x-b|^{pq}}
\ dm(x) \right)^{\frac{1}{q}}\,\left(\int\limits_{G}
Q^{q^{\prime}}(x)\ dm(x)\right)^{\frac{1}{q^{\prime}}}\,,
\end{equation}
где  $1/q+1/q^{\prime}=1$. Заметим, что первый интеграл в правой
части неравенства (\ref{eq13}) может быть подсчитан непосредственно.
Действительно, пусть для начала $q^{\prime}=\frac{n}{n-p}$ (и,
следовательно, $q=\frac{n}{p}.$) Ввиду теоремы Фубини будем иметь:
$$
\int\limits_{\varepsilon<|x-b|<\varepsilon_0}\frac{1}{|x-b|^{pq}} \
dm(x)=\omega_{n-1}\int\limits_{\varepsilon}^{\varepsilon_0}
\frac{dt}{t}=\omega_{n-1}\log\frac{\varepsilon_0}{\varepsilon}\,.
$$
В обозначениях леммы \ref{lem1} мы будем иметь, что при
$\varepsilon\rightarrow 0$
$$
\frac{1}{I^p(\varepsilon,
\varepsilon_0)}\int\limits_{\varepsilon<|x-b|<\varepsilon_0}\frac{Q(x)}{|x-b|^p}
\ dm(x)\leqslant \omega^{\frac{p}{n}}_{n-1}\Vert
Q\Vert_{L^{\frac{n}{n-p}}(G)}\left(\log\frac{\varepsilon_0}{\varepsilon}\right)
^{-p+\frac{p}{n}}\,\rightarrow 0\,,
$$
что влечёт выполнение соотношения (\ref{eq4!}).

\medskip
Пусть теперь $q^{\prime}>\frac{n}{n-p}$ (т.е.,
$q=\frac{q^{\prime}}{q^{\prime}-1}$). В этом случае
$$
\int\limits_{\varepsilon<|x-b|<\varepsilon_0}\frac{1}{|x-b|^{pq}} \
dm(x) = \omega_{n-1}\int\limits_{\varepsilon}^{\varepsilon_0}
t^{n-\frac{pq^{\prime}}{q^{\prime}-1}-1}dt\leqslant$$$$\leqslant
\omega_{n-1}\int\limits_{0}^{\varepsilon_0}
t^{n-\frac{pq^{\prime}}{q^{\prime}-1}-1}dt
=\frac{\omega_{n-1}}{n-\frac{pq^{\prime}}{q^{\prime}-1}}\varepsilon^{n-\frac{pq^{\prime}}{q^{\prime}-1}}_0,
$$
и, значит,
$$
\frac{1}{I^p(\varepsilon, \varepsilon_0)}
\int\limits_{\varepsilon<|x-b|<\varepsilon_0}\frac{Q(x)}{|x-b|^p} \
dm(x)\leqslant \Vert
Q\Vert_{L^{q^{\prime}}(G)}\left(\frac{\omega_{n-1}}{n-\frac{pq^{\prime}}{q^{\prime}-1}}
\varepsilon^{n-\frac{pq^{\prime}}{q^{\prime}-1}}_0\right)^{\frac{1}{q}}\left(\log\frac{\varepsilon_0}{\varepsilon}\right)^{-p}\,,
$$
что также влечёт выполнение соотношения (\ref{eq4!}). Таким образом,
необходимое утверждение вытекает из леммы  \ref{lem1}.
\end{proof}$\Box$

{\bf 5. Заключительные замечания} Легко построить пример
последовательности открытых дискретных отображений, не
удовлетворяющих соотношениям вида (\ref{eq2*A}), сходящихся к
отображению, не являющемуся нульмерным. Скажем, при $n=2$ достаточно
положить $f_m(z)=x+iy/m,$ $m\in {\Bbb N}.$ Заметим, что указанная
последовательность удовлетворяет соотношению
$$M(\Gamma)\leqslant m\cdot M(f(\Gamma))\,,$$
так как $K_O(x, f)=\frac{\Vert f^{\,\prime}(x)\Vert^2}{J(x, f)}=m,$
и, значит, постоянная $K$ в (\ref{eq2}) равна $m.$ Данная
последовательность сходится к отображению $f(z)=x,$ не являющемуся
нульмерным, причиной чего является отсутствие общей <<мажоранты>> в
(\ref{eq2}). В работе \cite{RSS} исследован вопрос о сходимости
гомеоморфизмов, удовлетворяющих <<обратным>> к (\ref{eq2*A})
неравенствам; здесь для гомеоморфности предельного отображения
достаточно только локальной интегрируемости функций $Q$ в
(\ref{eq2*A}). Можно ли условия 1)--3) в теореме \ref{th3} заменить
на условие локальной интегрируемости функции функции $Q,$ либо эти
условия являются точными в некотором смысле, нам неизвестно.

КОНТАКТНАЯ ИНФОРМАЦИЯ

\medskip
\noindent{{\bf Евгений Александрович Севостьянов} \\
Житомирский государственный университет им.\ И.~Франко\\
кафедра математического анализа, ул. Большая Бердичевская, 40 \\
г.~Житомир, Украина, 10 008 \\ тел. +38 066 959 50 34 (моб.),
e-mail: esevostyanov2009@mail.ru}

\end{document}